\documentclass[10pt,conference]{IEEEtran}

\def\real{{\mathchoice%
{\hbox{\rm\setbox1=\hbox{I}\copy1\kern-.45\wd1 R}}
{\hbox{\rm\setbox1=\hbox{I}\copy1\kern-.45\wd1 R}}
{\hbox{\scriptsize\rm\setbox1=\hbox{I}\copy1\kern-.45\wd1 R}}
{\hbox{\scriptsize\rm\setbox1=\hbox{I}\copy1\kern-.45\wd1 R}}}}
\def\Zint{{\mathchoice{\setbox1=\hbox{\sf Z}\copy1\kern-.75\wd1\box1}
{\setbox1=\hbox{\sf Z}\copy1\kern-.75\wd1\box1}
{\setbox1=\hbox{\scriptsize\sf Z}\copy1\kern-.75\wd1\box1}
{\setbox1=\hbox{\scriptsize\sf Z}\copy1\kern-.75\wd1\box1}}}
\newcommand{\complex}{ \hbox{\rm C\kern-0.45em\rule[.07em]{.02em}{.58em}%
\kern 0.43em}}

\newcommand{\be}{\begin{equation}}
\newcommand{\ee}{\end{equation}}
\newcommand{\beqr}{\begin{eqnarray}}
\newcommand{\eeqr}{\end{eqnarray}}
\newcommand{\beqrx}{\begin{eqnarray*}}
\newcommand{\eeqrx}{\end{eqnarray*}}
\newcommand{\ba}{\left[ \begin{array}}
\newcommand{\ea}{\\ \end{array} \right]}
\newcommand{\bi}{\begin{itemize}}
\newcommand{\ei}{\end{itemize}}

\newtheorem{theorem}{Theorem}

\newtheorem{algorithm}{Algorithm}

\usepackage{epsfig}
\usepackage{amsmath}
\usepackage{amssymb}

\def\w{{\bf w}}

\def\x{{\bf x}}
\def\y{{\bf y}}

\def\x{{\mathbf x}}

\def\w{{\bf w}}

\def\x{{\bf x}}
\def\y{{\bf y}}

\def\real{{\rm Re}\,}

\def\be{\begin{equation}}
\def\ee{\end{equation}}
\def\ba{\left[\begin{array}}
\def\ea{\end{array}\right]}

\def\w{{\bf w}}

\def\x{{\bf x}}
\def\y{{\bf y}}

\begin{document}
\title{Breaking through the Thresholds: an Analysis for Iterative Reweighted $\ell_1$ Minimization via the Grassmann Angle Framework }

\author{
\authorblockN{Weiyu Xu}
\authorblockA{Caltech EE \\
weiyu@systems.caltech.edu} \and
\authorblockN{M. Amin Khajehnejad}
\authorblockA{Caltech EE\\
amin@caltech.edu} \and
\authorblockN{A. Salman Avestimehr}
\authorblockA{Caltech CMI\\
avestim@caltech.edu} \and
\authorblockN{Babak Hassibi}
\authorblockA{ Caltech EE\\
hassibi@systems.caltech.edu}} \maketitle

\begin{abstract}

It is now well understood that $\ell_1$ minimization algorithm is able
to recover sparse signals from incomplete measurements
\cite{CT,DT,D} and  sharp recoverable sparsity thresholds have also
been obtained for the $\ell_1$ minimization algorithm. However, even
though iterative reweighted $\ell_1$ minimization algorithms or related
algorithms have been empirically observed to boost the recoverable
sparsity thresholds for certain types of signals, no rigorous
theoretical results have been established to prove this fact. In
this paper, we try to provide a theoretical foundation for analyzing
the iterative reweighted $\ell_1$ algorithms. In particular, we show
that for a nontrivial class of signals, the iterative reweighted
$\ell_1$ minimization can indeed deliver recoverable sparsity
thresholds larger than that given in \cite{DT,D}. Our results are
based on a high-dimensional geometrical analysis (Grassmann angle
analysis) of the null-space characterization for $\ell_1$ minimization
and weighted $\ell_1$ minimization algorithms.

\end{abstract}

\vspace{-0.15in}

\vspace*{0.25in} \noindent {\bf Index Terms:} compressed sensing,
basis pursuit, Grassmann angle, reweighted $\ell_1$ minimization,
random linear subspaces

\section{Introduction}
In this paper we are interested in compressed sensing problems.
Namely, we would like to find $\x$ such that
\begin{equation}
A\x=\y \label{eq:system}\vspace{-0.1in}
\end{equation}
where $A$ is an $m\times n$ ($m<n$) measurement matrix, $\y$ is a
$m\times 1$ measurement vector and $\x$ is an $n\times 1$ unknown
vector with only $k$ ($k<m$) nonzero components. We will further
assume that the number of the measurements is $m=\delta n$ and the
number of the nonzero components of $\x$ is $k=\zeta n$, where
$0<\zeta<1$ and $0<\delta<1$ are constants independent of $n$
(clearly, $\delta
>\zeta$).

A particular way of solving (\ref{eq:system}) which has recently
generated a large amount of research is called $\ell_1$-optimization
(basis pursuit) \cite{CT}. It proposes solving the following problem
\begin{eqnarray}
\mbox{min} & & \|\x\|_{1}\nonumber \\
\mbox{subject to} & & A\x=\y. \label{eq:l1}
\end{eqnarray}
Quite remarkably in \cite{CT} the authors were able to show that if
the number of the measurements is $m=\delta n$ and if the matrix $A$
satisfies a special property called the restricted isometry property
(RIP), then any unknown vector $\x$ with no more than $k=\zeta n$
(where $\zeta$ is an absolute constant which is a function of
$\delta$, but independent of $n$, and explicitly bounded in
\cite{CT}) non-zero elements can be recovered by solving
(\ref{eq:l1}). 
Instead of characterizing the $m\times n$ matrix $A$ through the RIP
condition, in \cite{DT,D} the authors assume that $A$ constitutes a
$k$-neighborly polytope. It turns out (as shown in \cite{DT}) that
this characterization of the matrix $A$ is in fact a necessary and
sufficient condition for (\ref{eq:l1}) to produce the solution of
(\ref{eq:system}). Furthermore, using the results of
\cite{VS}\cite{Schneider}\cite{henk}, it can be shown that if the
matrix $A$ has i.i.d. zero-mean Gaussian entries with overwhelming
probability it also constitutes a $k$-neighborly polytope. The
precise relation between $m$ and $k$ in order for this to happen is
characterized in \cite{DT} as well. 



In this paper we will be interested in providing the theoretical
guarantees for the emerging iterative reweighted $\ell_1$ algorithms
\cite{CWB07}. These algorithms iteratively updated weights for each
element of $\x$ in the objective function of $\ell_1$ minimization,
based on the decoding results from previous iterations. Experiments
showed that the iterative reweighted $\ell_1$ algorithms can greatly
enhance the recoverable sparsity threshold for certain types of
signals, for example, sparse signals with Gaussian entries. However,
no rigorous theoretical results have been provided for establishing
this phenomenon. To quote from \cite{CWB07}, ``any result
quantifying the improvement of the reweighted algorithm for special
classes of sparse or nearly sparse signals would be significant". In
this paper, we try to provide a theoretical foundation for analyzing
the iterative reweighted $\ell_1$ algorithms. In particular, we show
that for a nontrivial class of signals, (It is worth noting that
empirically, the iterative reweighted $\ell_1$ algorithms do not always
improve the recoverable sparsity thresholds, for example, they often
fail to improve the recoverable sparsity thresholds when the
non-zero elements of the signals are ``flat" \cite{CWB07}), a
modified iterative reweighted $\ell_1$ minimization algorithm can
indeed deliver recoverable sparsity thresholds larger than those
given in \cite{DT,D} for unweighted $\ell_1$ minimization algorithms.
Our results are based on a high-dimensional geometrical analysis
(Grassmann angle analysis) of the null-space characterization for
$\ell_1$ minimization and weighted $\ell_1$ minimization algorithms. The
main idea is to show that the preceding $\ell_1$ minimization
iterations can provide certain information about the support set of
the signals and this support set information can be properly taken
advantage of to perfectly recover the signals  even though the
sparsity of the signal $\x$ itself is large.

This paper is structured as follows. In Section \ref{sec:algorithm},
we present the iterative reweighted $\ell_1$ algorithm for analysis.
The signal model for $\x$  will be given in Section
\ref{sec:sigmodel}. In Section \ref{sec:claim} and Section
\ref{sec:probnull}, we will show how the iterative reweighted $\ell_1$
minimization algorithm can indeed improve recoverable sparsity
thresholds. Numerical results will be given in Section
\ref{sec:numerical}.

\section{The Modified Iterative Reweighted $\ell_1$ Minimization Algorithm}
\label{sec:algorithm}

Let $w_{i}^{t}$,  $i = 1, . . . , n$,  denote the weights for the
$i$-th element $\x_{i}$ of  $\x$ in the $t$-th iteration of the
iterative reweighted $\ell_1$ minimization algorithm and let
$\textbf{W}^{t}$ be the diagonal matrix with $w_{1}^{t}, w_{2}^{t} ,
... ,w_{n}^{t}$ on the diagonal. In the paper \cite{CWB07}, the
following iterative reweighted $\ell_1$ minimization algorithm is
presented:

\begin{algorithm} \cite{CWB07}

\begin{enumerate}
\item Set the iteration count $t$ to zero and $w_{i}^{t}=1$, $i=1, ..., n$.
\item Solve the weighted $\ell_1$ minimization problem
\begin{equation}
\x^{t} = \arg{ \min{ \|\textbf{W}^{t}\x\|}}~~\text{subject to}~~ \y
= A\x.
\end{equation}
\item Update the weights: for each $i = 1, . . . , n$,
\begin{equation}
w_{i}^{t+1}=\frac{1}{|\x_{i}^{t}|+\epsilon'},
\end{equation}
where $\epsilon'$ is a tunable positive number.

\item Terminate on convergence or when $t$ attains
a specified maximum number of iterations $t_\text{max}$. Otherwise,
increment $t$ and go to step 2.
\end{enumerate}
\label{alg:main}
\end{algorithm}

For the sake of tractable analysis, we will give another iterative
reweighted $\ell_1$ minimization algorithm , but it still captures the
essence of the reweighted $\ell_1$ algorithm presented in \cite{CWB07}.
In our modified algorithm, we only do two $\ell_1$ minimization
programming, namely we stop at the time index $t=1$.

\begin{algorithm}

\begin{enumerate}
\item Set the iteration count $t$ to zero and $w_{i}^{t}=1$, $i=1, ..., n$.
\item Solve the weighted $\ell_1$ minimization problem
\begin{equation}
\x^{t} = \arg{ \min{ \|\textbf{W}^{t}\x\|}}~~\text{subject to}~~ \y
= A\x.
\end{equation}
\item Update the weights:
find the index set $K' \subset \{1,2, ..., n\}$ which corresponds to
the largest $(1-\epsilon) \rho_{F}(\delta)\delta n$ elements of
$\x^{0}$ in amplitudes, where $0<\epsilon<1$ is a specified
parameter and $\rho_{F}(\delta)$ is the weak threshold for perfect
recovery defined in \cite{DT} using $\ell_1$ minimization  (thus
$\zeta= \rho_{F}(\delta)\delta$ is the weak sparsity threshold).
Then assign the weight $W_1=1$ to those $w_{i}^{t+1}$ corresponding
to the set $K'$ and assign the weight $W_2=W$, $W>1$, to those
$w_{i}^{t+1}$ corresponding to the complementary set
$\bar{K'}=\{1,2, ..., n\} \setminus K'$.

\item Terminate on convergence or when $t=1$. Otherwise, increment $t$ and go to step 2.
\end{enumerate}
\label{alg:modmain}
\end{algorithm}

This modified algorithm is certainly different from the algorithm
from \cite{CWB07}, but the important thing is that both algorithms
assign bigger weights to those elements of $\x$ which are more
likely to be 0.

\section{Signal Model for $\x$}
\label{sec:sigmodel}

In this paper, we consider the following model for the
$n$-dimensional sparse signal $\x$. First of all, we assume that
there exists a set $K \subset \{1,2,...,n\}$ with cardinality
$|K|=(1-\epsilon) \rho_{F}(\delta)\delta n$ such that each of the
elements of $\x$ over the set $K$ is large in amplitude. W.L.O.G.,
those elements are assumed to be all larger than $a_{1}>0$. For a
given signal $\x$, one might take such set $K$ to be the set
corresponding to the $(1-\epsilon) \rho_{F}(\delta)\delta n$ largest
elements of $\x$ in amplitude.

Secondly, (let $\bar{K}=\{1,2,...,n\} \setminus K$), we assume that
the $\ell_1$ norm of $\x$ over the set $\bar{K}$, denoted by
$\|\x_{\bar{K}}\|_1$, is upper-bounded by $\Delta$, though $\Delta$
is allowed to take a non-diminishing portion of the total $\ell_1$ norm
$\|\x\|_1$ as $n \rightarrow \infty$. We further denote the support
set of $\x$ as $K_{\text{total}}$ and its complement as
$\bar{K}_{\text{total}}$. The sparsity of the signal $\x$, namely
the total number of nonzero elements in the signal $\x$ is then
$|K_{\text{total}}|=k_{\text{total}}=\xi n$, where $\xi$ can be
above the weak sparsity threshold $\zeta= \rho_{F}(\delta)\delta$
achievable using the $\ell_1$ algorithm.

In the following sections, we will show that if certain conditions
on $a_1$, $\Delta$ and the measurement matrix $A$ are satisfied, we
will be able to recover perfectly the signal $\x$ using Algorithm
\ref{alg:modmain} even though its sparsity level is above the
sparsity threshold for $\ell_1$ minimization. Intuitively, this is
because the weighted $\ell_1$ minimization puts larger weights on the
signal elements which are more likely to be zero, and puts smaller
weights on the signal support set, thus promoting sparsity at the
right positions. In order to achieve this, we need some prior
information about the support set of $\x$, which can be obtained
from the decoding results in previous iterations. We will first
argue that the equal-weighted $\ell_1$ minimization of Algorithm
\ref{alg:modmain} can sometimes provide very good information about
the support set of signal $\x$.

\section{Estimating the Support Set from the $\ell_1$ Minimization }
\label{sec:claim}

Since the set $K'$ corresponds to the largest elements in the
decoding results of $\ell_1$ minimization, one might guess that most of
the elements in $K'$ are also in the support set $K_{total}$. The
goal of this section is to get an upper bound on the cardinality of
the set $\bar{K}_\text{total} \cap K'$, namely the number of zero
elements of $\x$ over the set $K'$ .  To this end, we will first
give the notion of ``weak" robustness for the $\ell_1$ minimization.

Let $K$ be fixed and $\x_K$, the value of
  $\x$ on this set, be also fixed. Then the solution produced by
  (\ref{eq:l1}), ${\hat \x}$, will be called {\em weakly robust} if,
  for some $C>1$ and all possible $\x_{{\bar K}}$, it holds that
\begin{equation*}
 \|(\x-\hat{\x})_{\bar{K}}\|_1
\leq \frac{2C}{C-1} \|\x_{\bar{K}}\|_1,
\end{equation*}
and
\begin{equation*}
 \|\x_K\|_1-\|\hat{\x}_K\|_{1} \leq
\frac{2}{C-1} \|\x_{\bar{K}}\|_1
\end{equation*}

The above ``weak'' notion of robustness allows us to bound the error
$\|\x-{\hat
  \x}\|_1$ in the following way. If the matrix $A_K$, obtained by
  retaining only those columns of $A$ that are indexed by $K$, has
  full column rank, then the quantity
\begin{equation*}
\kappa=\max_{A\w=0, \w \neq 0}
\frac{\|\w_K\|_1}{\|\w_{\bar{K}}\|_1},
\end{equation*}
must be finite ($\kappa<\infty$). In particular, since $\x-{\hat
\x}$ is in the null space of $A$ ($\y = A\x = A{\hat \x}$), we have
\begin{eqnarray*}
\|\x-{\hat \x}\|_1 & = & \|(\x-{\hat \x})_K\|_1+\|(\x-{\hat \x})_{\bar
  K}\|_1 \\
& \leq & (1+\kappa)\|(\x-{\hat \x})_{\bar K}\|_1 \\
& \leq &  \frac{2C(1+\kappa)}{C-1} \|\x_{\bar{K}}\|_1,
\end{eqnarray*}
thus bounding the recovery error. We can now give necessary and
sufficient conditions on the measurement matrix $A$ to satisfy the
notion of weak robustness for $\ell_1$ minimization.

\begin{theorem}

For a given $C>1$, support set $K$, and $\x_K$, the solution
$\hat{\x}$
  produced by (\ref{eq:l1}) will be {\em weakly robust} if, and only if,
$\forall\w\in\textbf{R}^n$ such that $A\w=0$, we have
 \begin{equation}
\|\x_K+\w_{K}\|_1+ \|\frac{\w_{\bar{K}}}{C}\|_1 \geq \|\x_K\|_1;
 \label{eq:wthmeq1}
 \end{equation}

\end{theorem}

\noindent \begin{proof}
Sufficiency: Let $\w=\hat{\x}-\x$, for which
$A\w=A(\hat{\x}-\x)=0$. Since ${\hat \x}$ is the minimum $\ell_1$ norm
solution, we have $\|\x\|_{1} \geq \|{\hat \x}\|_1 = \| \x+\w
\|_{1}$, and therefore $\|\x_K\|_{1} +\|\x_{\bar K}\|_{1} \geq
\|{\hat \x}_K\|_{1} +\|{\hat \x}_{\bar K}\|_{1}$. Thus,
\begin{eqnarray*}
\|\x_K\|_1-\|\x_K+\w_K\|_1 & \geq &
\|\w_{\bar{K}}+\x_{\bar{K}}\|_1-\|\x_{\bar{K}}\|_1 \\
& \geq & \|\w_{\bar{K}}\|_1-2\|\x_{\bar{K}}\|_1.
\end{eqnarray*}
But the condition (\ref{eq:wthmeq1}) guarantees that
\begin{equation*}
\|\w_{\bar{K}}\|_1\geq C(\|\x_K\|_1-\|\x_K+\w_K\|_1),
\end{equation*}
so we have
\begin{equation*}
 \|\w_{\bar{K}}\|_1
\leq \frac{2C}{C-1} \|\x_{\bar{K}}\|_1,
\end{equation*}
and
\begin{equation*}
\|\x_K\|_1-\|\hat{\x}_K\|_{1} \leq \frac{2}{C-1} \|\x_{\bar{K}}\|_1,
\end{equation*}
as desired.

Necessity: Since in the above proof of the sufficiency, equalities
can be achieved in the triangular inequalities, the condition
(\ref{eq:wthmeq1}) is also a necessary condition for the weak
robustness to hold for every $\x$. (Otherwise, for certain $\x$'s,
there will be $\x'=\x+\w$ with $\|\x'\|_1<\|\x\|_1$ while violating
the respective robustness definitions. Also, such $\x'$ can be the
solution to ($\ref{eq:l1}$)).
\end{proof}

We should remark (without proof for the interest of space) that for
any $\delta>0$, $1<\epsilon<1$, let $|K|=(1-\epsilon)
\rho_{F}(\delta)\delta n$, and suppose each element of the
measurement matrix $A$ is sampled from i.i.d. Gaussian distribution,
then there exists a constant $C>1$ (as a function of $\delta$ and
$\epsilon$), such that the condition (\ref{eq:wthmeq1}) is satisfied
with overwhelming probability as the problem dimension $n
\rightarrow \infty$. At the same time,the parameter $\kappa$ defined
above is upper-bounded by a finite constant (independent of the
problem dimension $n$) with overwhelming probability as $n
\rightarrow \infty$. These claims can be shown by using the
Grasamann angle approach for the balancedness property of random
linear subspaces in \cite{Allerton}. In the current version of our
paper, we would make no attempt to explicitly express the parameters
$C$ and $\kappa$.

In Algorithm \ref{alg:modmain}, after  equal-weighted $\ell_1$
minimization, we pick the set $K'$ corresponding to the
$(1-\epsilon) \rho_{F}(\delta)\delta$ largest elements in amplitudes
from  the decoding result $\hat{\x}$ (namely $\x^{0}$ in the
algorithm description) and assign the weights $W_{1}=1$ to the
corresponding elements in the next iteration of reweighted $\ell_1$
minimization. Now we can show that an overwhelming portion of the
set $K'$  are also in the support set $K_\text{total}$ of $\x$ if
the measurement matrix $A$ satisfies the specified weak robustness
property


\begin{theorem}
Supposed that we are given a signal vector $\x \in R^{n}$ satisfying
the signal model defined in Section \ref{sec:sigmodel}. Given
$\delta>0$,  and a measurement matrix $A$ which satisfies the weak
robustness condition in  (\ref{eq:wthmeq1})  with its corresponding
$C>1$ and $\kappa<\infty$, then the set $K'$ generated by the
equal-weighted $\ell_1$ minimization in Algorithm 2 contains at most
$\frac{2C}{(C-1)\frac{a_{1}}{2}} \|\x_{\bar{K}}\|_{1}+\frac{2C
\kappa}{(C-1)\frac{a_{1}}{2}} \|\x_{\bar{K}}\|_{1}$ indices which
are outside the support set of signal $\x$. \label{thm:key}
\end{theorem}

\noindent \begin{proof} Since the measurement matrix $A$ satisfies
the weak robustness condition for the set $K$ and the signal $\x$,

\begin{equation*}
 \|(\x-\hat{\x})_{\bar{K}}\|_1
\leq \frac{2C}{C-1} \|\x_{\bar{K}}\|_1.
\end{equation*}

By the definition of the $\kappa<\infty$, namely,
\begin{equation*}
\kappa=\max_{A\w=0, \w \neq 0}
\frac{\|\w_K\|_1}{\|\w_{\bar{K}}\|_1},
\end{equation*}

we have
\begin{eqnarray*}
 \|(\x-{\hat \x})_K\|_1
& \leq & \kappa\|(\x-{\hat \x})_{\bar K}\|_1.
\end{eqnarray*}

Then there are at most $\frac{2C}{(C-1)\frac{a_{1}}{2}}
\|\x_{\bar{K}}\|_{1}$ indices that are outside the support set of
$\x$ but have amplitudes larger than $\frac{a_{1}}{2}$ in the
corresponding positions of  the decoding result $\hat{\x}$ from the
equal-weighted $\ell_1$ minimization algorithm. This bound follows
easily from the facts that all such indices are in the set $\bar{K}$
and that $\|(\x-\hat{\x})_{\bar{K}}\|_1 \leq \frac{2C}{C-1}
\|\x_{\bar{K}}\|_1$.

Similarly, there are at most $\frac{2C \kappa}{(C-1)\frac{a_{1}}{2}}
\|\x_{\bar{K}}\|_{1}$ indices which are originally in the set $K$
but now have corresponding amplitudes smaller than $\frac{a_{1}}{2}$
in the decoded result $\hat{\x}$ of the equal-weighted $\ell_1$
algorithm.

Since the set $K'$ corresponds to the largest $(1-\epsilon)
\rho_{F}(\delta)\delta n$ elements of the signal $\x$, by combining
the previous two results, it is not hard to see that the number of
indices which are outside the support set of $\x$ but are in the set
$K'$ is no bigger than$\frac{2C}{(C-1)\frac{a_{1}}{2}}
\|\x_{\bar{K}}\|_{1}+\frac{2C \kappa}{(C-1)\frac{a_{1}}{2}}
\|\x_{\bar{K}}\|_{1}$.

\end{proof}

As we can see, Theorem \ref{thm:key} provides useful information
about the support set of the signal $\x$, which can be used in the
analysis for the weighted $l1$ minimization using the null-space
Grassmann Angle analysis approach for weighted $\ell_1$ minimization
algorithm \cite{isitweighted}.

\section{The Grassmann Angle Approach for the Reweighted $\ell_1$ Minimization}
\label{sec:probnull}
In the previous work \cite{isitweighted},  the authors have shown
that by exploiting certain prior information about the original
signal, it is possible to extend the threshold of sparsity factor
for successful recovery beyond the original bounds of \cite{DT,D}.
The authors proposed a nonuniform sparsity model in which the
entries of the vector $\x$ can be considered as $T$ different
classes, where in the $i$th class, each entry is (independently from
others) nonzero with probability $P_i$, and zero with probability
$1-P_i$. The signals generated based on this model will have around
$n_1P_1+\cdots + n_TP_T$ nonzero entries with high probability,
where $n_i$ is the size of the $i$th class. Examples of such signals
arise in many applications as  medical or natural imaging, satellite
imaging, DNA micro-arrays, network monitoring and so on. They prove
that provided such structural prior information is available about
the signal, a proper \emph{weighted} $\ell_1$-minimization strictly
outperforms the regular $\ell_1$-minimization in recovering signals
with some fixed average sparsity from under-determined linear i.i.d.
Gaussian measurements.

The detailed analysis in \cite{isitweighted} is only done for $T=2$,
and is based on the high dimensional geometrical interpretations of
the constrained weighted $\ell_1$-minimization problem:
\begin{equation*}
\min_{A\x = \y} \sum_{i=1}^n{w_i|\x_i|} 
\label{eq:weighetd l_1}
\end{equation*}
Let the two classes of entries be denoted by $K_1$ and $K_2$. Also,
due to the partial symmetry, for any suboptimal set of weights
$\{w_1,\cdots.\w_n\}$ we have the following

\begin{equation*}
\forall i\in\{1,2,\cdots,n\} ~~~w_i=\left\{\begin{array}{c}W_1
~if~i\in K_1\\ W_2 ~if~i\in K_2
\end{array}\right.
\label{eq:w_i's}
\end{equation*}

The following theorem is implicitly proven in \cite{isitweighted}
and more explicitly stated and proven in \cite{journalweighted}
\begin{theorem}
\label{thm:ISIT results} Let $\gamma_1=\frac{n_1}{n}$ and
$\gamma_2=\frac{n_2}{n}$. If $\gamma_1$, $\gamma_2$, $P_1$, $P_2$,
$W_1$ and $W_2$ are fixed, there exists a critical threshold
$\delta_c=\delta_c(\gamma_1,\gamma_1,P_1,P_2,\frac{W_2}{W_1})$,
totally computable, such that if $\delta=\frac{m}{n}\geq \delta_c$,
then a vector $\x$ generated randomly based on the described
nonuniformly sparse model can be recovered from the weighted
$\ell_1$-minimization of  \ref{eq:weighetd l_1} with probability
$1-o(e^{-cn})$ for some positive constant $c$.
\end{theorem}

In \cite{isitweighted} and \cite{journalweighted}, a way for
computing $\delta_c$ is presented which, in the uniform sparse case
(e.g $\gamma_2=0$) and equal weights, is consistent with the weak
threshold of Donoho and Tanner for almost sure recovery of sparse
signals with $\ell_1$-minimization.

In summary, given a certain $\delta$, the two different weights
$W_{1}$ and $W_{2}$ for weighted $\ell_1$ minimization,  the size of
the two weighted blocks, and also the number (or proportion) of
nonzero elements inside each weighted block, the framework from
\cite{isitweighted} can determine whether a uniform random
measurement matrix will be able to perfectly recover the original
signals with overwhelming probability. Using this framework we can
now begin to analyze the performance of the modified re-weighted
algorithm of section \ref{sec:algorithm}. Although we are not
directly given some prior information, as in the nonuniform sparse
model for instance, about the signal structure, one might hope to
infer such information after the first step of the modified
re-weighted algorithm. To this end, note that the immediate step in
the algorithm after the regular $\ell_1$-minimization is to choose the
largest $(1-\epsilon)\rho_F(\delta)\delta n$ entries in absolute
value. This is equivalent to splitting the index set of the vector
$\x$ to two classes $K'$ and $K''$, where $K'$ corresponds to the
larger entries. We now try to find a correspondence between this
setup and the setup of \cite{isitweighted} where sparsity factors on
the sets $K'$ and $\bar{K'}$  are known. We claim the following
upper bound on the number of nonzero entries of $\x$ with index on
$K'$
\begin{theorem}
There at at least $(1-\epsilon)\rho_F(\delta)\delta n -
\frac{4C(\kappa+1)\Delta}{(C-1)a_1}$ nonzero entries in $\x$ with
index on the set $K'$.
\end{theorem}

\begin{proof}
Directly from Theorem \ref{thm:key} and the fact that
$\|\x_{\bar{K}}\|_1\leq \Delta$.
\end{proof}

The above result simply gives us a lower bound on the sparsity
factor (ratio of nonzero elements) in the vector $\x_{K'}$
\begin{equation*}
P_1 \geq
1-\frac{4C(\kappa+1)}{(C-1)a_1\rho_F(\delta)\delta}\frac{\Delta}{n}
\label{eq:P_1}
\end{equation*}

Since we also know the original sparsity of the signal, $\|\x\|_0
\leq k_\text{total}$, we have the following lower bound on the
sparsity factor of the second block of the signal $\x_{\bar{K'}}$
\begin{equation*}
P_2 \leq \frac{k_\text{total}-(1-\epsilon)\rho_F(\delta)\delta n +
\frac{4C(\kappa+1)\Delta}{(C-1)a_1}}{n-(1-\epsilon)\rho_F(\delta)\delta
n}
\label{eq:P_2}
\end{equation*}

Note that if $a_1$ is large and $ 1 \gg \frac{\Delta}{a_{1}n} $
(Note however, we can let $\Delta$ take a non-diminishing portion
of $\|\x\|_{1}$, even though that portion can be very small),  then
$P_1$ is very close to 1. This means that the original signal is
much denser in the block $K'$ than in the second block $\bar{K'}$.
Therefore, as in the last step of the modified re-weighted
algorithm, we may assign a weight $W_1=1$ to all entries of $\x$ in
$K'$ and weight $W_2=W$, $W>1$ to the entries of $\x$ in $\bar{K'}$
and perform the weighted $\ell_1$-minimization. The theoretical results
of \cite{isitweighted}, namely Theorem \ref{thm:ISIT results}
guarantee that as long as $\delta >
\delta_c(\gamma_1,\gamma_2,P_1,P_2,\frac{W_2}{W_1})$ then the signal
will be recovered with overwhelming probability for large $n$. The
numerical examples in the next Section do show that the reweighted
$\ell_1$ algorithm can increase the recoverable sparsity threshold, i.e. $P_1\gamma_1 + P_2\gamma_2$.

\section{Numerical Computations on the Bounds}
\label{sec:numerical}
\begin{figure}
  \centering
  \includegraphics[width=0.45\textwidth]{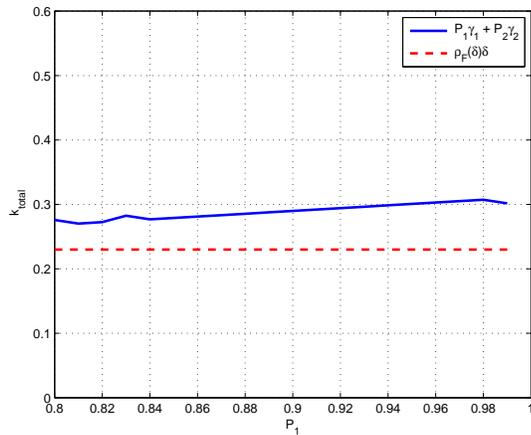}
   \caption{ \small Recoverable sparsity factor for $\delta=0.555$, when the modified re-weighted $\ell_1$-minimization algorithm is used.}
  \label{fig:sparsity curve}
\end{figure}

Using numerical evaluations similar to those in~\cite{isitweighted}, we demonstrate a strict improvement in the sparsity threshold from the weak bound of~\cite{DT}, for which our algorithm is guaranteed to succeed. Let $\delta=0.555$ and $\frac{W_2}{W_1}$ be fixed, which means that $\zeta=\rho_F(\delta)\delta$ is also given. We set $\epsilon=0.01$. The sizes of the two classes $K'$ and $\overline{K'}$ would then be $\gamma_1 n =(1-\epsilon)\zeta n$ and $\gamma_2 n = (1-\gamma_1)n$ respectively. The sparsity ratios $P_1$ and $P_2$ of course depend on other parameters of the original signal, as is given in equations (\ref{eq:P_1}) and (\ref{eq:P_2}). For values of $P_1$ close to $1$, we search over all pairs of $P_1$ and $P_2$ such that the critical threshold $\delta_c(\gamma_1,\gamma_2,P_1,P_2,\frac{W_2}{W_1})$ is strictly less than $\delta$. This essentially means that a non-uniform signal with sparsity factors $P_1$ and $P_2$ over the sets $K'$ and $\overline{K'}$ is highly probable to be recovered successfully via the weighted $\ell_1$-minimization with weights $W_1$ and $W_2$. For any such $P_1$ and $P_2$, the signal parameters ($\Delta$, $a_1$) can be adjusted accordingly. Eventually, we will be able to recover signals with average sparsity factor $P_1\gamma_1+P_2\gamma_2$ using this method. We simply plot this ratio as a function of $P_1$ in Figure \ref{fig:sparsity curve}. The straight line is the weak bound of~\cite{DT} for $\delta=0.555$ which is basically $\rho_F(\delta)\delta$.

\bibliographystyle{IEEEbib}

\end{document}